\documentclass[a4paper,tablecaptionabove,final]{article} %anstatt final auch draft sinnvoll scrartcl
\usepackage[utf8]{inputenc}
\usepackage[english]{babel}

\usepackage[usenames, dvipsnames]{color}
\usepackage{graphicx,enumerate} %enumerate, um \begin{enumerate}[a)]

\usepackage[intlimits]{amsmath}
\usepackage{amssymb, amsfonts}
\usepackage[amsmath,thmmarks,thref, hyperref]{ntheorem} %evtl nocht ,hyperref
\usepackage{dsfont} %\mathhds{R} zum richtigen Setzen der reellen Zahlen..

\usepackage[ngerman]{varioref} %zum besseren Verweis auf die Labels: \vref{},
\newcommand{\R}{\ensuremath{\mathds{R}}} %reelle Zahlen
\newcommand{\PP}{\ensuremath{\mathds{P}}} %P mit Doppelstrich

\newcommand{\E}{\ensuremath{\mathbf{E}\,}}

\theoremstyle{break}
\newtheorem{Theorem}{Theorem}[section]
\newtheorem{Satz}[Theorem]{Theorem}

\newtheorem{Definition}[Theorem]{Definition}

\theoremstyle{plain}

\newtheorem{Lemma}[Theorem]{Lemma}

\theoremstyle{nonumberplain} %nonumberbreak
\theoremseparator{:}
\theoremheaderfont{\normalfont\slshape}
\theorembodyfont{\normalfont}
\theoremsymbol{\ensuremath{\Box}}

\title{\Large{The Cameron-Martin Theorem for\\
 (p-)Slepian processes}}

\author{\normalsize \textbf{Wolfgang Bischoff\thanks{corresponding author:
Wolfgang Bischoff, Faculty of Mathematics and Geography,  D-85071 Eichst\"{a}tt, Germany;
E-Mail: wolfgang.bischoff@ku-eichstaett.de} \,  and Andreas Gegg
}\\[1ex]
\normalsize Faculty of Mathematics and Geography,\\
\normalsize Catholic University of Eichst\"{a}tt-Ingolstadt, \\
\normalsize }

\date{}

\graphicspath{ {Bilder/} }

\begin{document}

\maketitle

\begin{abstract}
  We show a Cameron-Martin theorem for Slepian processes $W_t:=\frac{1}{\sqrt{p}}(B_t-B_{t-p}), t\in [p,1]$, where $p\geq \frac{1}{2}$ and $B_s$ is Brownian motion. More exactly,  we determine the class of functions $F$ for which a density of $F(t)+W_t$ with respect to $W_t$ exists. Moreover, we prove an explicit formula for this density. p-Slepian processes are closely related to Slepian processes. p-Slepian processes play a prominent role among others in scan statistics and in testing for parameter constancy when data are taken from a moving window.
\end{abstract}

\section{Introduction} \label{cha:intro}
A Cameron-Martin theorem for a stochastic process is one of the most useful tools to solve problems the process is involved in. Let a stochastic process $X_{[a,b]}=(X_t)_{t\in [a,b]}$ with paths in $C[a,b]$, the set of real-valued and continuous functions on $[a,b]\subseteq \R$, and a deterministic function $F \in C[a,b]$ be given. Let $\PP^{X_{[a,b]}}$ and $\PP^{(X+F)_{[a,b]}}$ denote the distribution of $X_{[a,b]}$ and of $(X_t+F(t))_{t\in [a,b]}$ on $C[a,b]$, respectively. Then a Cameron-Martin theorem gives conditions on $F$ under which a density $\frac{d(\PP^{X_t+F(t)})}{d\PP^{X_t}}$ exists and, additionally, it gives an explicit formula for this density.

The first and best known result of this type is by \cite{c44}. They proved it for the standard Brownian motion $B_{[0,1]}:=(B_t)_{t\in [0,1]}$ with continuous paths. This Cameron-Martin theorem can be used, for example, to calculate optimal tests, see \cite{bm00}, and to estimate boundary crossing probabilities, see \cite{bh05}. Therefore, such a result of Cameron-Martin type is also of great interest for other stochastic processes. The following theorem, see Lifshits [12, Theorem 5.1], can be used as basis to get results of Cameron-Martin type for centered Gaussian processes. To this end, let $L^2(C[a,b],P)$ be the set of all real-valued and square-integrable functions on $C[a,b]$ with respect to a measure $P$ defined on the Borel-$\sigma$-Algebra of $C[a,b]$.

%%%%%%%%%%%%%%%%%%%%%%%%%%%%%%%%%%%%%%%%%%%%555
\begin{Satz}[An abstract Cameron-Martin theorem]  \label{th:GeneralCameronMartin}
Let $X_{[a,b]}:=(X_t)_{t\in[a,b]}$ be a centered Gaussian process with paths in $C[a,b]$, let $\mathcal{H}=\mathcal{H}_{X_{[a,b]}}\subseteq C[a,b]$ be the kernel of $\PP^{X_{[a,b]}}$ and $\|\cdot\|_{\mathcal{H}}$ its inherent norm. Then,
$$
\PP^{(X+h)_{[a,b]}} \text{ is absolutely continuous with respect to } \PP^{X_{[a,b]}} \iff h\in\mathcal{H}\; .$$
If $h\in\mathcal{H}$, then
$$\frac{d \PP^{(X+h)_{[a,b]}}}{d \PP^{X_{[a,b]}}}(g) =
\exp\left( - \frac12  \|h\|_{\mathcal{H}}^2 + z(g) \right)\quad \text{for }
\PP^{X_{[a,b]}}\text{-almost all } g \in C([a,b]),$$
where $z=z_{X_{[a,b]}}\in L^2(C[a,b],\PP^{X_{[a,b]}})$ is a linear functional fulfilling the equation
$$
\E\left(X_t \cdot z_{X_{[a,b]}}\right)=h(t)\;,\quad t\in[a,b].
$$
\end{Satz}
Hence, the kernel $\mathcal{H}=\mathcal{H}_{X_{[a,b]}}$ of the stochastic process $X_{[a,b]}$ together with its inherent norm $\|\cdot\|_\mathcal{H}$ and the functional $z_{X_{[a,b]}}$ must be determined to obtain an applicable result.

In the following, we use the notation stated in the next definition.

%%%%%%%%%%%%%%%%%%%%%%%%%%%%%%%%%%%%%%%%%%%%%%%%%%%%%%%%%%%%%%%%%%%%%%%
\begin{Definition}
Let $X_{[a,b]}:=(X_t)_{t\in[a,b]}$ be a centered Gaussian process with paths in $C[a,b]$ and let $\mathcal{H}=\mathcal{H}_{X_{[a,b]}}\subseteq C[a,b]$ be the kernel of $\PP^{X_{[a,b]}}$. Then, we say that a function
$$
z=z_{X_{[a,b]}}: \mathcal{H}_{X_{[a,b]}}\times C[a,b] \to \R
$$
fulfills condition A for $X_{[a,b]}$, if and only if for all $h\in \mathcal{H}_{X_{[a,b]}}$ the function $z(h,\cdot)\in L^2(C[a,b],\PP^{X_{[a,b]}}) $ is a linear functional fulfilling
\begin{align} \label{eq:z}
\E\left(X_t \cdot z(h,X_{[a,b]})\right)=h(t)\;,\quad t\in[a,b].
\end{align}
\end{Definition}
For the Brownian motion, we have the following result by \cite{c44}. To this end, let $\lambda$ denote the Lebesgue measure on $\R$ and let $L^2[a,b]:=L^2([a,b],\lambda)$ be the square-integrable functions on $[a,b]$ with respect to $\lambda$.

%%%%%%%%%%%%%%%%%%%%%%%%%%%%%%%%%%%%%%%%%%%%%%%%%%%%%%%%%%%%%%%%%%%%%%%%%%%%%%%%%%%%%%%%%%%%%%%%%%%%%%%%%%%%%%%%%%%%%%%%%%%
\begin{Satz} \label{BMCameronMartin}The kernel of $B_{[0,1]}$ is given by
$$
\mathcal{H}_{B_{[0,1]}} := \{ s_f\in C[0,1] \; | \; f\in L^2[0,1]\},
$$ \label{KernelHilbertBM}
where
$$
s_f(t):=\int_{[0,t]}f~d\lambda~, t\in [0,1],~\mbox{for $\lambda$-integrable } f:[0,1] \to \R.
$$
It is furnished with the norm
$$
\|s_f\|_{\mathcal{H}_{B_{[0,1]}}}:=\|f\|_{L^2[0,1]}, ~~f\in L^2[0,1].
$$
A function $z:=z_{B_{[0,1]}}$ fulfilling condition A for $B_{[0,1]}$ is given by the Wiener integral, see, for instance, \cite{mp10},
$$
z:\mathcal{H}_{B_{[0,1]}}\times C[0,1]\to\R,~~(s_f,b) \to \int_0^1 f\; db.
$$
\end{Satz}

Slepian and p-Slepian processes are defined and their close relations are discussed in the next section. In Sect. 3, results of Cameron-Martin type are established for p-Slepian processes. Several proofs are postponed to an appendix.

%%%%%%%%%%%%%%%%%%%%%%%%%%%%%%%%%%%%%%%%%%%%%%%%%%%%%%%%%%%%%%%%%%%%%%%%%%%%%%%
%%%%%%%%%%%%%%%%%%%%%%%%%%%%%%%%%%%%%%%%%%%%%%%%%%%%%%%%%%%%%%%%%%%%%%%%%%%%%%%
\section{Slepian and p-Slepian processes} \label{cha:intro}
The Slepian process X is the centered stationary Gaussian process with covariance function
 \begin{align*}
  C_X(s', s'+u)  = (1-u)^+\quad, ~0\leq s'\leq s'+u ,
 \end{align*}
where $t^+=\max(0,t), t\in\R.$ This process was introduced and studied in \cite{s61} and later in [16,17]. Afterward, it was handled in numerous theoretical and applied probabilistic models; see, e.g., [1-3,9-11,13,15].

Let us consider the Slepian process on time intervals $[1,b]$, where $1<b$ is any fixed constant. Let $\nabla_p~, p>0,$ be the backward difference operator with lag $p$, i.e., $\nabla_pF(t)=F(t)-F(t-p), t\in \R$, for functions $F: \R \to \R.$ Note that the Slepian process coincides in distribution with
$$
(\nabla_1B)_{[1,b]} := (\nabla_1B_t)_{t\in [1,b]}=(B_{t}-B_{t-1})_{t\in [1,b]}~,
$$
where $B_{t}$ is the standard Brownian motion with continuous paths in $[0,b]$. In Cressie [8, pp. 834], a slightly different expression of Slepian processes appeared in connection with scan statistics.
They were considered the processes
\begin{equation}  \label{eq:3}
\frac{1}{\sqrt{p}}(\nabla_pB)_{[p,1]}:=\frac{1}{\sqrt{p}}(\nabla_pB_t)_{t\in [p,1]}= \frac{1}{\sqrt{p}}\left( B_{t}-B_{t-p} \right)_{t\in[p,1]}
\end{equation}
for any fixed constant $p\in(0,1)$. We call this process $p$-Slepian process. Another application of p-Slepian processes is given by \cite{c95}. They consider moving sums of recursive residuals which are taken from windows of length $p$. By letting the number of residuals to $\infty$, they get the p-Slepian process $\frac{1}{\sqrt{p}}(\nabla_pB)_{[p,1]}$.

By a suitable scaling in time the Slepian process $(\nabla_1B)_{[1,b]}, b\in(1,\infty),$ can be transferred to the p-Slepian process $\frac{1}{\sqrt{p}}(\nabla_pB)_{[p,1]}, p\in(0,1)$, by putting $p=\frac{1}{b}.$
More exactly it holds $(\nabla_1B_u)_{u\in [1,b]}=\sqrt{b}
(\nabla_\frac{1}{b}B_\frac{u}{b})_{u\in [1,b]}$ in distribution.
%=\frac{1}{\sqrt{p}}(\nabla_pB_t)_{t\in [p,1]}$ with $p=\frac{1}{b}, t=\frac{u}{b}.$
The covariance function of a p-Slepian process is given by
 \begin{align*}
C_{\frac{1}{\sqrt{p}}(\nabla_pB)_{[p,1]}}(s,s+t) = \left(1-\frac{t}{p}\right)^+ \quad, p\leq s\leq s+t\leq 1.
 \end{align*}

%%%%%%%%%%%%%%%%%%%%%%%%%%%%%%%%%%%%%%%%%%%%%%%%%%%%%%%%%%%%%%%%%%%%%%%%%%%%%%%%%%%%%%%%%
%%%%%%%%%%%%%%%%%%%%%%%%%%%%%%%%%%%%%%%%%%%%%%%%%%%%%%%%%%%%%%%%%%%%%%%%%%%%%%%%%%%%%%%%%
%%%%%%%%%%%%%%%%%%%%%%%%%%%%%%%%%%%%%%%%%%%%%%%%%%%%%%%%%%%%%%%%%%%%%%%%%%%%%%%
\section{Cameron-Martin Theorem for p-Slepian processes} \label{sec:CMT for LF BM}
The function $\frac{1}{\sqrt{p}}\nabla_p, 0<p\leq 1,$ is linear. Hence, by \cite{l10} proposition 4.1, the kernel $
\mathcal{H}_{p-\text{Sl}} $ of the p-Slepian process $\frac{1}{\sqrt{p}}(\nabla_pB)_{[p,1]}$ is given by $\mathcal{H}_{p-Sl} = \frac{1}{\sqrt{p}}\nabla_p(\mathcal{H}_{B_{[0,1]}})= \{\frac{1}{\sqrt{p}}\nabla_p h| h\in \mathcal{H}_{B_{[0,1]}}\}.$ By some calculations we get
\begin{align}
\mathcal{H}_{p-\text{Sl}} \;
&= \; \left\{ \frac{1}{\sqrt{p}}\nabla_p s_f:[p,1]\to \R \; | \; f\in L^2[0,1] \right\}  \nonumber \\
& = \left\{c + s_g:[p,1]\to \R \; | \; c\in\R \; \text{, } \;  g\in L^2[p, 1]
\right\}.\nonumber
\end{align}
More exactly, we have for $f\in L^2[0,1]$
\begin{equation}
\frac{1}{\sqrt{p}}\nabla_p s_f(t)=\frac{1}{\sqrt{p}}s_f(p)+s_{\frac{1}{\sqrt{p}}\nabla_p f}(t), t\in [p,1].
\end{equation}
It is much more complicated to determine the inherent norm of $\mathcal{H}_{p-\text{Sl}} $ since $\nabla_p$ is not injective.  For the following, the information of $\nabla_p s_{f}$ is important. For $f,g\in L^2[0,1]$, it holds true
\begin{align}
&\forall t\in [p,1]:\nabla_p s_{f}(t)=s_f(p)+\int_p^tf(s)-f(s-p) ds=\nabla_p s_{g}(t)
 \Longleftrightarrow\nonumber\\
&f(t)-f(t-p)=g(t)-g(t-p) ~\lambda-a.s. \mbox{ for all }t\in [p,1],~~ s_f(p)=s_g(p)\nonumber.
\end{align}
We prove the following result in the appendix.
%%%%%%%%%%%%%%%%%%%%%%%%%%%%%%%%%%%%%%%%%
\begin{Lemma} \label{lemma norm}
Let $1/2\leq p \leq 1$. Then the kernel $\mathcal{H}_{p-\text{Sl}} $ is furnished with its inherent norm
\begin{align}
 & \left\|\frac{1}{\sqrt{p}}\nabla_p s_{f}\right\|_{p-Sl}^2 = \frac{1}{p}\inf_{g\in L^2[0,1]:\nabla_p s_{f}=\nabla_p s_{g}}\|g\|^2_{L^2[0,1]}\\
&=\frac{1}{2p(3p-1)}(2s_{f}(p)+\delta(1-p))^2+\frac{1}{2p}\|(f(t)-f(t-p)) _{t\in [p,1]}\|^2_{L^2[p,1]}, \nonumber
   \end{align}
where $\delta=\frac{1}{1-p}(s_f(1)-s_f(p)-s_f(1-p))$.
This minimum is attained at the function
\begin{align}
f^*(t)=&\left(\frac{1}{3p-1}(s_{f}(p)+\frac{1-p}{2}\delta)+\frac{1}{2}(-f(t+p)+f(t))\right){\bf 1}_{[0,1-p]}(t)\nonumber\\
&+\frac{1}{2p-1}(s_f(p)-\frac{1-p}{3p-1}(s_{f}(p)-(2p-1)\delta)){\bf 1}_{(1-p,p)}(t)\nonumber\\
&+\left(\frac{1}{3p-1}(s_{f}(p)+\frac{1-p}{2}\delta)+\frac{1}{2}(f(t+p)-f(t))\right){\bf 1}_{[p,1]}(t), t\in[0,1].\nonumber
  \end{align}
\end{Lemma}
%%%%%%%%%%%%%%%%%%%%%%%%%%%%%%%%%5
Finally, we need the function  $z:=z_{\frac{1}{\sqrt{p}}(\nabla_pB)_{[p,1]}}$ fulfilling condition A for the p-Slepian process to be in the position to state a result of Cameron-Martin type.
\begin{Lemma} \label{lemma fct z}
Let $1/2\leq p\leq 1$. The function $z:=z_{\frac{1}{\sqrt{p}}(\nabla_pB)_{[p,1]}}$ defined by
\begin{align}\label{Z WSP}
z:&\mathcal{H}_{\frac{1}{\sqrt{p}}(\nabla_pB)_{[0,1]}}\times C[p,1]\to\R,~~\left(\frac{1}{\sqrt{p}}s_f(p)+s_{\frac{1}{\sqrt{p}}\nabla_p f},\frac{1}{\sqrt{p}}\nabla_pb\right) \mapsto\\
&\frac{3s_f(p)+s_{\nabla_p f}(1)}{2(3p-1)}\left(\nabla_pb(p)+\nabla_pb(1)\right) + \frac{1}{2}\int_p^{1} \nabla_p f(s) \; d (\nabla_pb(s))
\end{align}
fulfills condition A for $\frac{1}{\sqrt{p}}(\nabla_pB)_{[p,1]}$. Note that the integral
$$
\int_p^{1} \nabla_p f(s) \; d (\nabla_pB(s))=\int_p^{1} \nabla_p f(s) \; dB(s)-\int_p^{1} \nabla_p f(s) \; dB(s-p)
$$
has to be understood as Wiener Integral.
\end{Lemma}
The proof is given in the appendix. By the above lemmas, we obtain the following result.

%%%%%%%%%%%%%%%%%%%%%%%%%%%%%%%%%%%%%%%%%%%%555
\begin{Satz}[Cameron-Martin theorem for p-Slepian processes]  \label{th:GeneralCameronMartin}
Let $1/2\leq p\leq 1$. It holds true
\begin{align*}
&\PP^{(\frac{1}{\sqrt{p}}\nabla_pB+h)_{[p,1]}} \text{ is absolutely continuous with respect to } \PP^{(\frac{1}{\sqrt{p}}\nabla_pB)_{[p,1]}}\\
 &\quad \quad \quad \iff h\in\mathcal{H}_{p-\text{Sl}}\; .
\end{align*}
For $h=\frac{1}{\sqrt{p}}s_f(p)+s_{\frac{1}{\sqrt{p}}\nabla_p f}(\cdot)\in\mathcal{H}_{p-\text{Sl}}$ it holds true for $\PP^{(\frac{1}{\sqrt{p}}\nabla_pB)_{[p,1]}}$-almost all $\frac{1}{\sqrt{p}}\nabla_pb \in C[p,1]$:
\begin{align*}
&\frac{d \PP^{(\frac{1}{\sqrt{p}}\nabla_pB+h)_{[p,1]}}}{d \PP^{(\frac{1}{\sqrt{p}}\nabla_pB)_{[p,1]}}}\left(\frac{1}{\sqrt{p}}\nabla_pb\right) \\
&=\exp\left( - \frac{1}{4p}\left(\frac{(2s_{f}(p)+\delta(1-p))^2}{3p-1} +\|(f(\cdot)-f(\cdot-p))\|^2_{L^2[p,1]}\right)\right.\\
 &~~~~\left.+ \frac{3s_f(p)+s_{\nabla_p f}(1)}{2(3p-1)} (\nabla_pb(p)+\nabla_pb(1)) + \frac{1}{2}\int_p^{1} \nabla_p f(s) \; d (\nabla_pb(s)) \right),
\end{align*}
where $\delta=\frac{1}{1-p}(s_f(1)-s_f(p)-s_f(1-p))$.
\end{Satz}
%%%%%%%%%%%%%%%%%%%%%%%%%%%%%%%%%%%%%%%%%%%%%%%%%%%%%%%%%%%%%%%%%%%%%%%%%%%%%%%
%%%%%%%%%%%%%%%%%%%%%%%%%%%%%%%%%%%%%%%%%%%%%%%%%%%%%%%%%%%%%%%%%%%%%%%%%%%%%%%
%%%%%%%%%%%%%%%%%%%%%%%%%%%%%%%%%%%%%%%%%%%%%%%%%%%%%%%%%%%%%%%%%%%%%%%%%%%%%%%
\section{Appendix} \label{sec:app}

%%%%%%%%%%%%%%%%%%%%%%%%%%%%%%%%%%%%%%%%%%%%%%%%%%%%%%%%%%%%%%%%%%%%%%%%%%%%%%%%%%%%%%%%
{\it Proof of Lemma \ref{lemma norm}}
Let $\frac{1}{\sqrt{p}}\nabla_p s_{f}\in \mathcal{H}_{p-\text{Sl}}, f\in L^2[0,1]$, be arbitrary and let $1/2\leq p\leq 1$. There are several places where the case $p=1/2$ must be considered separately. We do not do this since the specifications for $p=1/2$ are simpler than the following considerations for $1/2<p\leq 1$. By proposition 4.1 of \cite{l10}, we get
\begin{align}
\|\frac{1}{\sqrt{p}}\nabla_p s_{f}\|_{p-Sl} =& \frac{1}{\sqrt{p}}\inf_{g\in L^2[0,1]:\nabla_p s_{f}=\nabla_p s_{g}}\|s_g\|_{\mathcal{H}_{B_{[0,1]}}}\nonumber\\
=& \frac{1}{\sqrt{p}}\inf_{g\in L^2[0,1]:\nabla_p s_{f}=\nabla_p s_{g}}\|g\|_{L^2[0,1]},\nonumber
\end{align}
where the second equality follows by Theorem \ref{BMCameronMartin}.
Each function $f\in L^2[0,1]$ can be written $\lambda-a.s.$ in the form
\begin{align}\label{form}
f(t)=\alpha{\bf 1}_{[0,1-p]}(t)+a(t)+\beta{\bf 1}_{(1-p,p)}(t)+b(t)+\gamma{\bf 1}_{[p,1]}(t)+c(t), t\in [0,1],
\end{align}
where
\begin{align*}
&\alpha=\frac{1}{1-p}s_f(1-p), a(t)={\bf 1}_{[0,1-p]}(t)(f(t)-\alpha),\\
&\beta=\frac{1}{2p-1}(s_f(p)-s_f(1-p)), b(t)={\bf 1}_{(1-p,p)}(t)(f(t)-\beta),\\
&\gamma=\frac{1}{1-p}(s_f(1)-s_f(p)), c(t)={\bf 1}_{[p,1]}(t)(f(t)-\gamma).
\end{align*}
Note that the six summands in (\ref{form}) build orthogonal functions in $L^2[0,1]$. In the following, we use the consequence of this representation at several places without citing it explicitly. Furthermore, we have
\begin{align}
s_f(p)=\alpha (1-p)+ \beta (2p-1), f(t)-f(t-p)=\gamma-\alpha+c(t)-a(t), t\in [p,1].\nonumber
 \end{align}
%%%%%%%%%%%%%%%%%%%%%%%%%%%%%%%%%%%%%%%%%%%%%%%%%%%%%%
Next, we consider $g\in L^2[0,1]$ of the specific form
\begin{align}
g(t)=\alpha {\bf 1}_{[0,1-p]}(t)+\beta{\bf 1}_{(1-p,p)}(t)+(\alpha+\delta){\bf 1}_{[p,1]}(t), t\in[0,1], \alpha,\beta,\delta \in \R.\nonumber
 \end{align}\label{stepfct}
Hence,
$$
s_g(p)=\alpha (1-p)+\beta (2p-1) \Leftrightarrow \beta=\frac{1}{2p-1}(s_g(p)-\alpha (1-p))
$$
and
$$
g(t)-g(t-p)=\delta, t\in[p,1].
$$
Thus, $\nabla_p(g)(t)=s_g(p)+\delta t, t\in [p,1]$. Let $g_0(t)=\alpha_0 {\bf 1}_{[0,1-p]}(t)+\beta_0{\bf 1}_{(1-p,p)}(t)+(\alpha_0+\delta_0){\bf 1}_{[p,1]}(t)\in L^2[0,1]$ be fixed and let $g\in L^2[0,1]$ with $\nabla_p(g)(t)=\nabla_p(g_0)(t)=s_{g_0}(p)+\delta_0 t, t\in [p,1]$. Then,
$$
g(t)=\alpha {\bf 1}_{[0,1-p]}(t)+\frac{1}{2p-1}(s_{g_0}(p)-\alpha (1-p)){\bf 1}_{(1-p,p)}(t)+(\alpha+\delta_0){\bf 1}_{[p,1]}(t)
$$
and the square of its norm is given by
 \begin{align}
  \|g\|^2= \frac{1-p}{2p-1} (\alpha^2 (3p-1)-2\alpha (s_{g_0}(p)-(2p-1)\delta_0))+\frac{s_{g_0}(p)^2}{2p-1}+\delta_0^2(1-p).\nonumber
  \end{align}
This norm is minimal if and only if
 \begin{align}
\alpha=\frac{1}{3p-1}(s_{g_0}(p)-(2p-1)\delta_0). \nonumber
\end{align}
Hence, we obtain after some calculation
 \begin{align}
  &\|\nabla_p s_{g_0}\|^2_{p-Sl}= \frac{1}{3p-1}(2s_{g_0}(p)^2+2(1-p)s_{g_0}(p)\delta_0+\delta_0^2(1-p)p).\nonumber
  \end{align}
Next, we consider a general function
$$
f=\alpha{\bf 1}_{[0,1-p]}(t)+a(t)+\beta{\bf 1}_{(1-p,p)}(t)+b(t)+(\alpha+\delta){\bf 1}_{[p,1]}(t)+c(t)\in L^2[0,1].
$$
By the above considerations we obtain
\begin{align*}
  \|\nabla_p s_f\|_{p-Sl}^2
  =& \frac{1}{3p-1}(2s_{f}(p)^2+2(1-p)s_{f}(p)\delta+\delta^2(1-p)p)\\
  &+\min (\|a(t)\|^2+\|b(t)\|^2+\|c(t)\|^2)\nonumber
   \end{align*}
where the minimum is taken over all $a:[0,1-p]\to\R, b:(1-p,p)\to\R, c:[p,1]\to\R $ with $s_a(1-p)=s_b(2p-1)=s_c(1)=0$ and $c(t)-a(t-p)=f(t)-f(t-p)-\delta$. It holds true for $t\in[p,1]$ fixed
\begin{align}
\min|a(t-p)|^2+|c(t)|^2%\nonumber\\
   &=|\frac{1}{2}(-c(t)+a(t-p))|^2+|\frac{1}{2}(c(t)-a(t-p))|^2\nonumber\\
   &=\frac{1}{2}|(c(t)-a(t-p))|^2%\nonumber\\
    =\frac{1}{2}|f(t)-f(t-p)-\delta|^2\nonumber,
\end{align}
where the minimum is taken over all $a(t-p), c(t)\in\R$ with $c(t)-a(t-p)=f(t)-f(t-p)-\delta$.  Hence,
\begin{align}
  \|\nabla_p s_f\|_{p-Sl}^2 \nonumber
    = &\frac{1}{3p-1}(2s_{f}(p)^2+2(1-p)s_{f}(p)\delta+\delta^2(1-p)p)\nonumber\\
     &+\min \|a(t)\|_{L^2[0,1]}^2+\|c(t)\|_{L^2[0,1]}^2.\nonumber\\
    =&\frac{1}{3p-1}(2s_{f}(p)^2+2(1-p)s_{f}(p)\delta+\delta^2(1-p)p)\nonumber\\
    &+\frac{1}{2}\|(f(t)-f(t-p)-\delta)_{t\in[p,1]}\|_{L^2[p,1]}^2\nonumber\\
    =&\frac{1}{3p-1}(2s_{f}(p)^2+2(1-p)s_{f}(p)\delta+\delta^2(1-p)p)\nonumber\\
    &+\frac{1}{2}(\|(f(t)-f(t-p))_{t\in [p,1]}\|_{L^2[p,1]}^2
     -(1-p)\delta^2)\nonumber\\
    =&\frac{1}{2(3p-1)}(2s_{f}(p)+\delta(1-p))^2 \nonumber\\
    &+\frac{1}{2}\|(f(t)-f(t-p))_{t\in[p,1]}\|_{L^2[p,1]}^2.\nonumber
   \end{align}
This minimum is obtained at the function
\begin{align}
f^*(t)=&\left(\frac{1}{3p-1}(s_{f}(p)+\frac{1-p}{2}\delta)+\frac{1}{2}(-f(t+p)+f(t))\right){\bf 1}_{[0,1-p]}(t)\nonumber\\
&+\frac{1}{2p-1}(s_f(p)-\frac{1-p}{3p-1}(s_{f}(p)-(2p-1)\delta)){\bf 1}_{(1-p,p)}(t)\nonumber\\
&+\left(\frac{1}{3p-1}(s_{f}(p)+\frac{1-p}{2}\delta)+\frac{1}{2}(f(t+p)-f(t))\right){\bf 1}_{[p,1]}(t), t\in[0,1].\nonumber
  \end{align}

%%%%%%%%%%%%%%%%%%%%%%%%%%%%%%%%%%%%%%%%%%%%%%%%%%%%%%%%%%%

{\it Proof of Lemma \ref{lemma fct z}}
We have to prove Eq. (\ref{eq:z}) for the p-Slepian-process. To this end, let $p\leq t\leq 1, f\in L^2[0,1]$ and $(\nabla_p s_f)(\cdot)=s_f(p)+s_{\nabla_p f}(\cdot)\in \mathcal{H}_{p-\text{Sl}}.$ It holds true for $t\in[p,1]$:
\begin{align*}
&\E\left(\frac{1}{\sqrt{p}}(\nabla_pB)_{[p,1]}(t)\cdot \frac{3s_f(p)+s_{\nabla_p f}(1)}{2(3p-1)} \left((\nabla_pB)_{[p,1]}(p)+(\nabla_pB)_{[p,1]}(1)\right) \right)\\ &\quad\quad\quad+\E\left(\frac{1}{\sqrt{p}}(\nabla_pB)_{[p,1]}(t)\cdot \frac{1}{2}\int_p^{1} \nabla_p f(s) \; d(\nabla_pB)_{[p,1]}(s))\right) \\
&=\frac{3s_f(p)+s_{\nabla_p f}(1)}{2\sqrt{p}(3p-1)} \E\left((B_{[0,1]}(t)-B_{[0,1]}(t-p))(B_{[0,1]}(p)+B_{[0,1]}(1)-B_{[0,1]}(1-p)) \right)\\
&\quad\quad\quad+\frac{1}{2\sqrt{p}}\E\left((B_{[0,1]}(t)-B_{[0,1]}(t-p))\int_p^{1} \nabla_p f(s) \; d (B_{[0,1]}(s)-B_{[0,1]}(s-p))\right) \\
&=\frac{3s_f(p)+s_{\nabla_p f}(1)}{2\sqrt{p}}-\frac{1}{2\sqrt{p}}\left(s_{\nabla_p f}(1)-s_{\nabla_p f}(p)\right)= \frac{1}{\sqrt{p}}s_f(p)+s_{\frac{1}{\sqrt{p}}\nabla_p f}(t).
\end{align*}

\end{document}